\documentclass[11pt]{article}
\usepackage{amsmath, amssymb, amsthm, verbatim,enumerate,bbm,tikz}
\usepackage{indentfirst}

\usepackage{graphicx}

\definecolor{myblue}{RGB}{80,80,160}
\definecolor{mygreen}{RGB}{80,160,80}

\date{\today}
\allowdisplaybreaks

\theoremstyle{plain}
\newtheorem{theorem}{Theorem}[section]
\newtheorem{lemma}[theorem]{Lemma}

\newtheorem{corollary}[theorem]{Corollary}

\newcommand{\Var}{{\bf Var \,}}
\newcommand{\BN}{\mathbb N }
\newcommand{\BP}{{\bf P}}
\newcommand{\BE}{{\bf E}}
\newcommand{\CL}{{\cal L}}
\newcommand{\CE}{{\cal E}}
\newcommand{\BI}{{ \bf I}}

\newcommand{\tendD}{ \overset {d}{\longrightarrow} }

\title{Random matrices: Law of the iterated logartihm}

\author{ Asaf Ferber
\thanks{Department of Mathematics, Yale University. Email:
asaf.ferber@yale.edu.} \and Daniel Montealegre \thanks{Department of Mathematics, Yale University. Email: daniel.montealegre@yale.edu}\and Van Vu \thanks{Department of
Mathematics, Yale University. Email: van.vu@yale.edu.} }

\begin{document}
\maketitle

\begin{abstract}
The theory of random matrices contains many  central limit theorems. We have central limit theorems for
eigenvalues statistics, for the log-determinant and log-permanent,   for limiting distribution of individual eigenvalues in the bulk, and many others.

In this notes, we discuss the following problem: Is it possible to prove the law of the iterated logarithm? We illustrate this possibility by showing that this is indeed the case for the log of the permanent of random Bernoulli matrices and pose
open questions concerning several other matrix parameters.
\end{abstract}

\section{Introduction}
Let $\xi_i, i=1,2 \dots $  be an infinite sequence of iid random variables with mean 0 and variance 1. The most important result  in probability, the  classical central limit theorem
(CLT), asserts that for $Y_n:= \sum_{i=1}^n \xi_i$, one has

$$ \frac{ Y_n  }{\sigma_n  }  \overset {d}{\longrightarrow} N(0,1),$$ where $\sigma_n := \sqrt {\Var Y_n }  =\sqrt n $, $N(0,1)$ denotes the standard gaussian and
$\tendD$ denotes convergence in distribution.


In  the 1920s, Khinchin and Kolmogorov proved the famous law of the  iterated logarithm (LIL)   \cite{Khinchin, Kolmogorov} , which asserts that

$$ \BP \left[ \limsup_{n\rightarrow \infty} \frac{ Y_n  }{\sigma_n \sqrt {2 \log \log n }} =1 \right] =1 . $$

What is remarkable about the LIL  is that it takes into account the correlation between the $\xi_i$, which was not needed for the CLT. For instance, the CLT holds if one considers
a triangular array  $\xi_{ij},  j \le i $ of iid variables and define
$Y_n =\sum_{j=1}^n \xi_{nj} $.  All logarithms  in this paper have natural base.


In this paper, we consider the possibility of proving LIL in random matrix setting.  Let  $$ M_{\infty} := ( \xi_{ij})_{ij} , $$  where $\xi_{ij}$ are independent random variables, be an infinite matrix and
 $M_n$ be its  principle minor formed by the first $n$ rows and columns (having entries $\xi_{ij} , 1\le i, j \le n$). Let $Y_n$ be a parameter of $M_n$ which satisfies the central limit theorem, that is,  there are real numbers $\mu_n, \sigma_n$ such that

$$ \frac{ Y_n -\mu_n }{ \sigma_n}  \overset {d}{\longrightarrow}  N(0,1). $$

\noindent  A natural question is  whether $Y_n$ satisfies LIL, namely

\begin{equation}  \label{LIL1}   \BP  \left[ \limsup   \frac{ Y_n -\mu_n }{  \sigma_n \sqrt { 2 \log \log n } }  =1\right]  =1. \end{equation}  

\noindent The above setting is non-hermitian.  For hermitian (or symmetric) matrices, one naturally adds the condition
that $\xi_{ij}  =\bar \xi_{ji} $ (or $\xi_{ij} = \xi_{ji} $) and an appropriate condition on the diagonal entries.

\vskip2mm

The main result  of this paper is the  LIL  for the permanent of random matrices.

\begin{theorem}  \label{theorem:main1}
Consider the infinite (non-hermitian) matrix $M_{\infty}$ with entries $\xi_{ij}$ being iid Bernoulli variables (taking values 0 and 1 with probability $1/2$). Let  $X_n$ be the permanent of $M_n$ and $Y_n = \log  X_n$. Then  $Y_n$ satisfies the LIL,

\begin{equation}  \label{LIL1}   \BP \left[ \limsup   \frac{ Y_n -  \log \frac{n!}{2^n}  }{ \sqrt {2 \log \log n } }  =1\right]  =1. \end{equation}
\end{theorem}

The corresponding CLT was proved earlier by Janson \cite[Theorem 14]{Janson} and
also by Rempala et al. \cite{Rem}  in a more general form.

Apparently, Theorem \ref{theorem:main1} is only the tip of an iceberg. To motivate further investigation in this direction,   let us state a few concrete  open problems.

\vskip2mm

\noindent {\bf LIL for log-determinant.}   Let $\xi_{ij}$ be iid sub-gaussian random variables with mean 0 and variance 1. Let $Y_n  = \log |\det M_n |$. Nguyen et al. \cite{NV}  proved that

$$ \frac{ Y_n  -  \frac{1}{2} \log (n-1) !} {\sqrt {2 \log n } } \tendD N(0,1). $$

\noindent  {\bf Problem.}  Does  $Y_n$ satisfy the LIL ?

\vskip2mm

\noindent  {\bf LIL for linear statistic of eigenvalues.} Consider the Hermitian model with the upper diagonal entries
$\xi_{ij}, i < j$ be iid sub-gaussian random variables with mean 0 and variance 1, and the  diagonal entries be iid sub-gaussian random variables with mean 0 and variance 2. Let $\phi$ be a
{\it nice} test function, and  define $Y_n:= \sum_{i=1}^n \phi (\lambda_i ) $. It is well known that  $Y_n$ satisfies the CLT.
There is a large literature on this phenomenon (with many different definitions of {\it nice}); see, for instance,  \cite[Section 18.4]{Pastur} for details.

\noindent {\bf Problem.}  Does  $Y_n$ satisfy the LIL ?

\vskip2mm

We note that  for some parameters, it could happen that the right normalization is not $\sqrt { \log \log n }$ (and it is natural to view the above  questions  in this broader sense). In \cite{PZ} Paquette  et al. considered  the infinite GUE matrix   and defined
$Y_n := (\lambda_n -  2 \sqrt n ) n^{1/6} $, where $\lambda_n$ is the largest eigenvalue.   They  showed, for an explicit constant $c$, that

\begin{equation}  \label{LIL2}   \BP \left[ \limsup   \frac{ Y_n }{  \log^{2/3} n  }  =c \right]  =1. \end{equation}

Thus, one obtained a fractional logarithm, rather than iterated logarithm, law.  On the other hand, this particular  $Y_n$ does not satisfy the CLT, either.  Another relevant result is \cite[Proposition 5.4]{Bourgade} which
studied  the LIL in a very different setting.

\vskip2mm

The rest of the paper is organized as follows.
In Section \ref{tool}, we state a slightly more general version of our theorem and
the main lemma behind its proof.  Section \ref{lemma} is devoted to the verification of the lemma.
In the last section, Section \ref{theorem},  we prove the theorem.  In order to main the flow of the arguments, we delay the proofs of several technical
estimates  to the appendix.

\section {A more general statement and the main technical lemma}  \label{tool}

\subsection{ Bernoulli matrices with arbitrary density}

Our theorem still holds if we allow the Bernoulli random variables to have mean $p$, for any constant $0 < p < 1$.
Moreover, it also has a combinatorial interpretation. To see this, one needs to define  a  bipartite graph $G_n$ associated with $M_n$. Technically,
$G_n$ is the bipartite graph with vertices indexed by the rows and columns of $M_n$ and a vertex $i$ (in the ``row" color class) is
connected to a vertex $j$ (in the ``column" class) if and only if the corresponding matrix entry $\xi_{ij} =1$.  The permanent of $M_n$ is precisely the number of perfect matchings in
$G_n$.

The random matrix $M_{\infty}$ (with mean $p$)  then gives rise to an infinite bipartite graph $G(\BN, \BN, p )$, where the two color classes consist of natural numbers, and
any two vertices $i,j$ (from different classes) are connected  independently with probability $p$.  Let $G(n,n,p)$  be the finite graph spanned by the first $n$ vertices in each color class.
This way, the  general version of Theorem \ref{theorem:main1} can be combinatorially formulated  as follows

\begin{theorem}\label{main thm 2}Let $0 < p < 1$ be a constant and $X_n$ be the  number of perfect matchings in  $G(n,n,p)$. Set $Y_n := \log X_n$. Then
\begin{align}\label{eqtn main thm 2}
\BP  \left[\lim\sup_{n\rightarrow \infty}\frac{ Y_n- \log (n! p^n )   }{ \sqrt{2\log \log n}  \sqrt { \frac{1-p}{p}} }=1\right]=1
\end{align}
\end{theorem}

This setting is more convenient for our proof, which relies on   combinatorial estimates. In particular, our main tool will be the concentration result discussed in the next section.

\subsection{  Concentration of the number of perfect matchings}

Let $m$ be a natural number. We define  a new model of random bipartite graph, denoted by  $G(n,n,m)$,  as follows. Consider two color classes with  $n$ vertices each (labelled by numbers from $1$ to $n$ as usual).
The edges of $G(n,n,m)$ are a (uniformly) random subset of exactly $m$ elements of the set of all possible $n^2$ edges between the two color classes.
Let $X_{n,m}$ denote the number of perfect matchings in $G(n,n,m)$.

\begin{lemma}\label{concentration 2}  Let $0 < \delta < 1/2$ be a constant. There is a constant $C$ depending on $\delta$ such that  for
 any $ \delta n^2 \le m \le (1-\delta)n^2 $, and $k = o(n^{1/3} )$, we have

\begin{equation} \label{conceq}   \BE X_{n,m}^k  \leq  C^k ( \BE  X_{n,m} ) ^k . \end{equation}

\end{lemma}

For more information  about random graphs and matchings, we refer the reader to \cite{JLR}.  Using Lemma \ref{concentration 2}, Markov's bound implies that for all $K \ge C$

$$\BP(   X_{n,m} \ge  K  \BE X_{n,m} ) \le (C/K)^k . $$

 By taking $\delta:=  \min \{ p/2, (1-p)/2 \}$, $k = 4\log n$ and $K=  Ce$,  we obtain the following corollary

\begin{corollary} \label{deviationbound}  Let $0 < p< 1$ be a constant.
There is a constant $K$ (depending on $p$) such that  for
 any $ \frac{p}{2}  n^2 \le m \le \frac{1+p}{2} n^2    $

$$  \BP  ( X_{n.m}  \ge  K \BE X_{n,m} ) \leq  n^{-4}. $$

\end{corollary}

\section{Proof of the concentration lemma}\label{lemma}

We denote by $K_{n,n}$ the complete bipartite graph (on the vertex set of $G(n,n,m)$) and let $\mathcal P$ to denote the set of all perfect matchings in $K_{n,n}$. Clearly, we have 

$$| \mathcal P |  = n! . $$

For each $P\in \mathcal P$, let $X_P$ to denote the indicator random variable for the event ``$P$ appears in $ G(n,n,m)$".
 It is easy to see that

  \begin{equation} \label{formula1}  \BE X_P =\frac{(m)_n}{(n^2)_n},  \end{equation}  where $(N)_n := N(N-1) \dots (N-n+1) $.
  Thus,

$$ \BE X_{n,m} = n!  \frac{ (m)_n} { (n^2) _n } . $$

\noindent A routine calculation (see the Appendix)  shows that

\begin{equation} \label{Xexpectation}
\frac{ (m)_n} {(n^2)_n} = p_m^n  \exp  \left( - \frac{1 -p_m}{p_m } +O(1/n) \right),
\end{equation}

where $p_m:=\frac{m}{n^2}$.  In general, for any fixed bipartite graph
  $H$ with $h$ edges, the probability that  $G(n,n,m)$ contains $H$ is precisely

$$   \frac{(m)_{h}}{{(n^2)}_{h}}.$$

\noindent We will make a repeated use of the following estimate which its simple proof appears in the Appendix
\begin{equation} \label{estimate1}  (N)_{\ell} = N^{\ell} \exp\left(- \frac {\ell (\ell-1) }{2N }+ o(1) \right) \end{equation} for all $N, \ell$ such that $ \ell =o(N^{2/3} )$.

 Thinking of $H$ as the (simple) graph formed by the union of perfect matchings $P_1, \dots, P_k$, observing that $X_H=X_{P_1}\cdots X_{P_k}$, we obtain that

\begin{equation}   \label{rearranged}
  \BE X_{n,m}^k = \sum_{P_1,...,P_k\in \mathcal P} \BE [X_{P_1} \dots X_{P_k}] = \sum_{a=0}^{(k-1)n}M(a)\frac{(m)_{kn-a}}{{(n^2)}_{kn-a}},
\end{equation}
where $M(a)$ is the number of (ordered) $k$-tuples
$(P_1,...,P_k)\in \mathcal P^k$, whose union contains exactly $kn-a$ edges. Our main task is to bound $M(a)$ from above.

 Fix $a$ and let $\CL:=\CL(a)$ be the set of all  sequences  $L:= \ell_2,\ldots,\ell_k$   of non-negative integers where  $$ \ell_2 + \dots + \ell_k = a. $$   For each sequence $L = \ell_2, \dots, \ell_k$,  let  $N_{L} $ be the number of $k$-tuples  $(P_1,\ldots,P_k)$ such that
 for every $2 \le t \le k$,  we have  $|P_t \cap (\cup_{j<t}P_j)|=\ell_t$. Clearly, we have

 $$M(a)  =\sum _{ L \in \CL } N_L . $$

 We construct a $k$-tuple in $N_L$ according to the following algorithm:

 \begin{itemize}
\item Let $P_1$ be an arbitrary perfect matching.

\item Suppose that $P_1,\ldots,P_{t-1}$ are given, our aim is to  construct $P_t$. Pick $\ell_t$ edges to be in $P_t \cap \cup_{j=1}^{t-1}P_j$ as follows: first, pick a subset $B_{1,t}$ of $\ell_t$ vertices from the first color class (say $V_1$). Next,
 from each vertex pick an edge which appears in
$\cup_{j=1}^{t-1}P_j$ so that the chosen edges form a matching. Let us denote the obtained partial matching by $E_t$, and observe that $|E_t|=\ell_t$, and that $B_{2,t}:=\left(\cup E_t\right)\cap V_2$ is a set of size $\ell_t$ (where $V_2$ denotes the second color class).

\item  Find a perfect matching $M_t$ between $V_1 \backslash B_{1,t} $ and $V_2 \backslash B_{2,t} $ which has an empty intersection with $\cup_{j=1}^{t-1} P_j$, and set $P_t:= E_t \cup M_t $.
\end{itemize}
\begin{center}
\includegraphics[scale=.5]{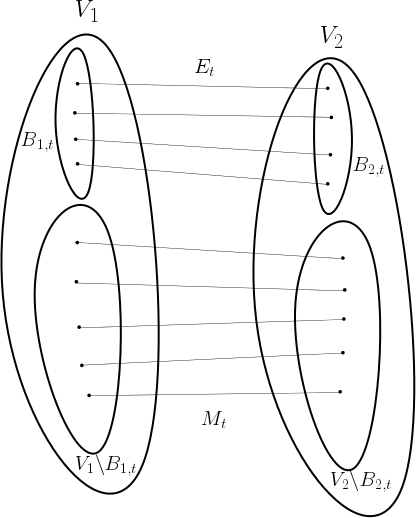}
\end{center}
Next, we wish to analyze the algorithm. There are $n!$ ways to choose $P_1$. Having chosen $P_1,\ldots,P_{t-1}$, there are ${n \choose {\ell_t} }$ ways to choose $B_{1,t}$.  Each vertex in $B_{1,t}$ has at most $t-1$ different edges in $\cup_{j=1}^{t-1}P_j$. Thus, the number of ways to choose
$E_t$ is at most $(t-1)^{\ell_t} $.  Moreover, once $B_{1,t}$ and $B_{2,t} $ are defined, the number of ways to choose $M_t$ is at most  $(n-\ell_t) !$. This way, we obtain

$$N_L \le n !  \prod_ {t=2}^{k}  {n \choose {\ell_t} } (t-1)^{\ell_t}  (n-\ell_t)! =  n !  \prod_{t=2}^k n!  \frac{(t-1)^{\ell_t} } {\ell_t ! } = (n !)^k  \prod_{t=2}^k   \frac{(t-1)^{\ell_t} } {\ell_t ! }. $$

\noindent By the multinomial identity  and the definition of the set $\CL$,

$$\sum_{L \in \CL}   \prod_{t=2}^k   \frac{(t-1)^{\ell_t} } {\ell_t ! }  =  \frac{1}{a!} (1+ \dots + (k-1)) ^a = \frac{ { k \choose 2}^a}{a!}. $$

\noindent Therefore

\begin{equation} \label{bound1} M(a) = \sum_{L \in \CL} N_L \le   (n!)^k \sum_{ L \in \CL } \prod_{t=2}^k   \frac{(t-1)^{\ell_t} } {\ell_t ! }  = (n!)^k  \frac{ {k \choose 2} ^a}{a!}. \end{equation}

This estimate is sufficient in the case $a$ is relatively large. However, it is too generous in the case $a$ is small (the main contribution in LHS of \eqref{rearranged}  comes from this case). In order to sharpen the bound, we refine the estimate on the number of possible $M_t$'s that one can choose in the last step of the algorithm, call this number $\mathcal M_t$ (clearly, $\mathcal M_t$ also depends on the $B_{i,t}$s and we estimate a worse case scenario).
Let $G_t$ be the bipartite graph between $V_1 \backslash B_{1,t} $ and $V_2 \backslash B_{2,t}$ formed by the edges which are not in $\cup_{j=1}^{t-1}P_j$. For each $v \in V_1 \backslash B_{1,t} $,  let $d_v$ be its degree in $G_t$.  By
Bregman-Minc inequality (see the Appendix)

$$\mathcal M_t \le \prod_{ v \in V_1 \backslash B_{1,t} }  (d_v !)^{1/d_v} . $$

\noindent It is clear from the definition that for each $v$

$$d  := n- \ell_t -(t-1) \le d_v \le n -\ell_t : = D$$

Call a vertex $v$ {\it good} if $d_v = d $ and {\it bad} otherwise. It is easy to see that $v$ is good iff it has exactly $t-1$ different edges in $\cup_{j=1}^{t-1} P_j$ and none of these edges hits  $B_{2,t} $. It follows that the number of good vertices is at least

$$ n- \ell_t (t-1) - \sum_{j=2}^{t-1} \ell_j  \ge n -a(k-1) -a = n - ka . $$

Since $(d!)^{1/d} $ is monotone increasing, it follows that

$$\mathcal M_t \le  (d !) ^{\frac{n-ka}{d} }   (D! )^{ \frac{ ka -\ell_t }{D } }. $$

Comparing to the previous bound of $(n-\ell_t) ! $, we gain a factor of

\begin{equation} \label{bound2}   \frac { (d !) ^{\frac{n-ka}{d} }   (D! )^{ \frac{ ka -\ell_t }{D } } }{ (n-\ell_t) !  } = \left[  \frac{ (d!)^{1/d} } { (D!)^{1/D } } \right] ^{n-ka}. \end{equation}

\noindent A routine calculation (see Appendix) shows that whenever $ka= o( n)$, the RHS is

\begin{equation}  \label{bound3}  (1 + o(1)) e^{- (t-1) }. \end{equation}

\noindent Thus, for such values of $a$, we have

\begin{equation} \label{bound4} M(a)  \le  (n!)^k  \frac{ {k \choose 2} ^a}{a!}  \prod_{t=2}^k (1+o(1)) e^{-(t-1)}  < 2^k \exp\left(- \frac{k(k-1)}{2}  \right)   (n!)^k  \frac{ {k \choose 2} ^a}{a!} ,  \end{equation}  where the constant 2 can be replaced by any constant larger than 1.

Now  we are ready to bound  $\BE X_{n,m} ^k$. Recall (\ref{rearranged})

$$ \BE X_{n,m}^k =\sum_{a=0}^{(k-1)n}M(a)\frac{(m)_{kn-a}}{{(n^2)}_{kn-a}} . $$

\noindent We split the RHS as

$$ \sum_{a=0}^{T}M(a)\frac{(m)_{kn-a}}{{(n^2)}_{kn-a}}+\sum_{a=T+1 }^{(k-1)n}M(a)\frac{(m)_{kn-a}}{{(n^2)}_{kn-a}} = S_1 +S_2. $$

\noindent where $T= pek^2$.  The assumption $k^3 =o(n)$ of the lemma guarantees that $kT =o(n)$. Let  $p_m := \frac{m}{n^2} $. By \eqref{bound4} and \eqref{estimate1} and a routine calculation, we have

$$ S_1= \sum_{a=0}^{T}M(a)\frac{(m)_{kn-a}}{(n^2)_{kn-a}}
\leq \frac{2^k(n!)^kp_m^{nk}}{e^{{k\choose 2}}}\exp\left(-\frac{k^2(1-p_m)}{2p_m}+o(1)\right)\sum_{a=0}^{T}\frac{({k\choose 2})^{a}}{a!}p_m^{-a}. $$

\noindent On the other hand,

$$ \sum_{a=0}^{T}\frac{({k\choose 2})^{a}}{a!}p_m^{-a} <  \sum_{a=0}^{\infty }\frac{({k\choose 2})^{a}}{a!}p_m^{-a}  = e^{ { k \choose2 } /p_m }, $$ so

$$S_1 \le  \frac{2^k  (n!)^kp_m ^{nk}}{e^{{k\choose 2}}}\exp\left(-\frac{k^2(1-p_m)}{2p_m}+o(1)\right)e^{{k\choose 2}/p_m}=
  C_1^k(n!)^k p_m^{nk}, $$

 \noindent where $C_1$ is a constant depending on $p$. (In fact we can replace the constant $2$ by any constant larger than 1 in the definition of $C_1$; see the remark following \eqref{bound4}).
  To bound $S_2$, we use \eqref{bound1} and \eqref{estimate1} to obtain

 $$S_2 = \sum_{a >T}  M(a)\frac{(m)_{kn-a}}{(n^2)_{kn-a}}
\leq (n!)^kp_m^{nk}  \exp\left(-\frac{k^2(1-p_m)}{2p_m}+o(1)\right)\sum_{a >T } \frac{({k\choose 2})^{a}}{a!}p_m^{-a}. $$

Notice that we no longer have the term $\frac{2^k}{ e^{{ k \choose 2}}}$. However, as $a$ is large, there is a much better way to bound
$\sum_{a >T } \frac{({k\choose 2})^{a}}{a!}p_m^{-a}. $ Stirling's approximation yields

$$ \sum_{a >T } \frac{({k\choose 2})^{a}}{a!}p_m^{-a} \le \sum_{a >T} \left(\frac{ek^2}{ 2  p_m a } \right)^a < \sum_{a >T} \left(\frac{1}{2} \right)^2 = o(1). $$

\noindent It follows that

$$S_2= o( (n!)^k p_m^{nk}) , $$ and thus is negligible for our needs.  Therefore,

$$\BE  X_{n,m} ^k= S_1 + S_2 \le C_1^k (n!) p_m^{nk}. $$

\noindent Finally, note that \eqref{Xexpectation} implies

$$(\BE X_{n,m})^k = (n!)^k p_m^{nk} \exp \left(\frac{k(p_m-1  )}{p_m } +O(k/n)\right)  \le C_2^k   (n!)^k p_m^{nk}, $$ where
for $C_2= \exp( \frac{p_m-1 }{p_m } + o(1) )$.    We conclude the proof of \eqref{conceq} by setting $C:= C_1 C_2^{-1} $.

\section{Proof of Theorem \ref{main thm 2} }\label{theorem}

\subsection{Upper Bound}\label{upper bound section}

We need to prove that for  any fixed  $\varepsilon>0$

\begin{equation} \label{upper}
\BP  \left[\frac{\log X_n-\log  (n! p^n) }{ \sqrt { \frac{1-p}{p} }}\geq (1+\varepsilon)\sqrt{2 \log \log n^2}\text{ for infinite many }n\right]=0. \end{equation}

We combine Corollary \ref{deviationbound} with an argument from \cite{Janson}.
By Corollary \ref{deviationbound}, there is a constant $K$ such that  for all $ \frac{p}{2} n^2 \le  m  \le \frac{1+p}{2} n^2 $

$$   X_{n,m}  \le K \BE X_{n,m}  $$ with probability at least $1 -n^{-4} $. Taking $\log$,  we conclude  that with the same probability

\begin{equation} \label{upper1}
 Y_{n,m} \le \log \BE  X_{n,m} + \log K. \end{equation}

\noindent Recalling the estimate \eqref{Xexpectation}, we have

$$\log \BE X_{n,m}  =\log (n! p_m^n) +  \frac{p_m-1}{2p_m } + o(1), $$ with $p_m := \frac{m}{n^2}$.

\noindent The  RHS can be written as

$$ \log(n!)+n\log\frac{m}{n^2}-\frac{n^2}{2}\left(\frac{1}{m}-\frac{1}{n^2}\right)+ o(1). $$

Let $E_n$ be the random variable that counts the number of edges in $G(n,n,p)$. By  conditioning on $E_n=m$ and using the union bound
(over the range $\frac{p}{2}n^2 \le m \le \frac{1+p}{2} n^2$), we can conclude that  with probability at least  $1 -n^{-2} $

$$ \BI_{\CE}  \log X_n \le  \BI_{\CE} \Big( \log(n!)+n\log\frac{E_n}{n^2}-\frac{n^2}{2}\left(\frac{1}{E_n}-\frac{1}{n^2} \right) + \log K +  o(1) \Big), $$

\noindent where $X_n$ denotes  the number of perfect matchings in $G(n,n,p)$,  and $\BI_{\CE}$ is the indicator of the event $\CE$ that
$G(n,n,p)$ has at least $\frac{p}{2} n^2$ and at most $\frac{1+ p}{2} n^2 $ edges. By Chernoff's bound, $\BI_{\CE}= 1$ with probability at least $1 -n^{-2} $. By the union bound

\begin{align} \label{upper0}  \log X_n \le \Big( \log(n!)+n\log\frac{E_n}{n^2}-\frac{n^2}{2}\left(\frac{1}{E_n}-\frac{1}{n^2}\right) + O(1) \Big),
\end{align}  with probability at least $1 - 2n^{-2} $.

Note that with  probability at least $1- n^{-2} $,  $E_n= pn^2 + O (n\log^2 n)$, in which case   $\frac{n^2}{2E_n} =O(1)$.
Again by the union bound, we have with probability at least $1- 3 n^{-2} $,

\begin{equation} \label{upper1}  \log X_n \le \Big( \log(n!)+n\log\frac{E_n}{n^2} + O(1) \Big).
\end{equation}

 Let  $E_n^*:=(E_n- \BE [E_n])/\sqrt{\Var (E_n)}$ , we have
\begin{align}
\log \frac{E_n}{n^2}&=\log \left(\frac{\sqrt{\Var(E_n)}E_n^*}{n^2}+\frac{\BE E_n}{n^2}\right)\nonumber\\
&=\log\left(\left(\frac{p(1-p)}{n^2}\right)^{1/2}E_n^*+p\right)\nonumber\\
&=\log\left(p\left(\frac{1-p}{p}\right)^{1/2}\frac{E_n^*}{n}+p\right)\nonumber\\
&=\log p+\log \left(1+\left(\frac{1-p}{p}\right)^{1/2}\frac{E_n^*}{n}\right)\nonumber\\
&=\log p+\left(\frac{1-p}{p}\right)^{1/2}\frac{E_n^*}{n}+O(1/n^2).  \nonumber
\end{align}

\noindent Plugging the last estimate into  \eqref{upper1} we obtain,  with probability at least  $1 -3n^{-2} $
\
\[
\log X_n \le \log (n! p^n) +\left(\frac{1-p}{p}\right)^{1/2}E_n^* + O(1), \nonumber
\]  or equivalently
\begin{align}\label{comparison}
\frac{\log X_n-\log (n! p^n) }{\sqrt {\frac{1-p}{p}}} =E_n^*+ O(1).
\end{align}

Since $\sum _n n^{-2} < \infty$, we have, by Borell-Cantelli lemma (see the Appendix) that the event in \eqref{comparison} holds with probability 1 for all sufficiently large $n$.
On the other hand, by Kolmogorov-Khinchin theorem,  $E_n^*$ satisfies LIL and thus
\[
E_n^*\leq (1+\varepsilon/2)\sqrt{2\log\log n^2} \le   (1+ 2\varepsilon/3 )\sqrt{2\log\log n}
\]happens with probability $1$ for all sufficiently large  $n$.  Note that $E_n$ is the sum of $n^2$ iid random variables, and thus we have $\log \log n^2$ here instead of
$\log \log n$.  Finally, for all sufficiently large $n$, $ (\varepsilon/3)\sqrt{2\log\log n }$ is larger than the error term
$O(1)$, and we have
\[
\frac{\log X_n-\log (n! p^n) }{\sqrt {\frac{1-p}{p}}}   \leq (1+\varepsilon)\sqrt{2\log\log n },
\] proving the upper bound.

\subsection{Proof of the Lower bound}

For the lower bound we need to show that there exists a sequence $n_k, k=1, 2 \dots$ of indices such that with probability 1,

\[
\frac{\log X_{n_k} -\log (n_k ! p^{n_k} ) }{\sqrt {\frac{1-p}{p}}}   \ge (1- \varepsilon)\sqrt{2\log\log n_k },
\]

\noindent holds for infinitely many $k$.

By Theorem \cite[Theorem 14]{Janson}, we know that $\Var X_{n,m} = O\left( \frac{1}{n} (\BE X_{n,m})^2\right)$.  Markov's bound then implies that

\begin{equation} \label{lower tail}  X_{n,m} \ge  \frac{3}{2} \BE X_{n,m} \end{equation}  with probability at least $1 -O(1/n)$.  This bound
 is sufficient here, as we only need to consider a very sparse subsequence.
From  the standard proof of LIL for sum of iid random variables \cite{Khinchin, Kolmogorov}, we see that there is a sequence $\{n_k\}:=\{c^k\}$ (where $c$ is an integer larger than 1) for which we have:
\[
E_{n_k}^*\geq (1-\varepsilon/2)\sqrt{2\log \log n_k^2}\geq (1-\varepsilon/2)\sqrt{2\log \log n_n}
\]
happens infinitely often with probability one. Restricting ourselves to this subsequence and repeating the calculation in the previous section, we obtain
for every $n_k$

$$ \label{comparison2}
\frac{\log X_{n_k} -\log (n_k ! p^{n_k} )  }{ \sqrt { \frac{1-p}{p}}} \ge E_{n_k}^*+ O(1), $$  with probability
$1- O( n_k ^{-1} ) $.

\noindent Let $A_k$ denote the event that equation \eqref{lower tail} fails for $n_k$. Then
\[
\BP A_k=O(1/c^k)
\]so in particular we have
\[
\sum_{k} \BP A_k <\infty.
\] By  Borel-Cantelli lemma (see Appendix), we have with probability equal to $1$ that for infinitely many $k$ the following two estimates holds.

\begin{itemize}
\item $E_{n_k}^*\geq (1-\varepsilon/2)\sqrt{2\log \log n_k}$.
\item$\frac{\log X_{n_k} -\log (n_k ! p^{n_k} ) }{ \sqrt { \frac{1-p}{p}} }=E_{n_k}^*+O(1). $
\end{itemize}

The lower bound now follows as for every large enough $k$, as  $(\varepsilon/2)\sqrt{2\log \log n_k}$ is greater than the error term $O(1) $.

\section{Appendix}
\textbf{Approximation lower factorial:} Let $N, \ell$ be such that $\ell=o(N^{2/3})$. Then,
\begin{align}
(N)_\ell&=N(N-1)\cdots (N-\ell+1) \nonumber\\
&=N^{\ell}\prod_{i=0}^{\ell-1}(1-i/N)\nonumber\\
&=N^{\ell}\prod_{i=0}^{\ell-1}e^{-i/N+O(i^2/N^2)}\nonumber\\
&=N^{\ell}\exp\left(\sum_{i=0}^{\ell-1}-i/N+O(i^2/N^2)\right)\nonumber\\
&=N^{\ell}\exp\left(-\frac{\ell(\ell-1)}{2N}+O(\ell^3/N^2)\right)\nonumber\\
&=N^{\ell}\exp\left(-\frac{\ell(\ell-1)}{2N}+o(1)\right)\nonumber
\end{align}as claimed.
\\
\\
\textbf{Computation of equation \eqref{bound2}}: We are going to use the following upper and lower bounds for the factorial:
\[
\sqrt{2\pi s}(s/e)^s\leq s!\leq \sqrt{2\pi s}(s/e)^se^{1/12s}
\]Hence,
\begin{align}
\left[  \frac{ (d!)^{1/d} } { (D!)^{1/D } } \right] ^{n-ka}&\leq\left[\frac{(\sqrt{2\pi d}(d/e)^de^{1/12d})^{1/d}}{(\sqrt{2\pi D}(D/e)^D)^{1/D}}\right]^{n-ka}\nonumber\\
&=\left[(1+O(n^{-2}))\frac{d(2\pi d)^{1/2d}}{D(2\pi D)^{1/2D}}\right]^{n-ka}\nonumber\\
&=(1+O(n^{-1})\left[\frac{(2\pi d)^{1/2d}}{(2\pi D)^{1/2D}}\right]^{n-ka}\left[\frac{d}{D}\right]^{n-ka}\nonumber\\
&=(1+o(1))\left[1-\frac{t-1}{n-\ell_t}\right]^{n-ka}\nonumber\\
&=(1+o(1))e^{t-1}\nonumber	
\end{align}as desired.  (Here we use the assumption that  $ka =o(n)$.)

\begin{lemma}[Borel-Cantelli Lemma]
Let $(A_i)_{i=1}^{\infty}$ be a sequence of events. Then
\begin{enumerate}[$(a)$]
\item If $\sum_{k=1}^{\infty} \BP \left[A_k\right]<\infty$, then
$$\BP \left[A_k \textrm{ holds for infinitely many }k\right]=0.$$
\item If $\sum_{k=1}^{\infty}  \BP \left[A_k\right]=\infty$ and in addition all
the $A_k$'s are independent, then
$$\BP \left[A_k \textrm{ holds for infinitely many }k\right]=1.$$
\end{enumerate}
\end{lemma}

\begin{theorem}[Bregman-Minc inequality; \cite{Bregman}]\label{bregman}Let $G$ be a bipartite graph with two color classes $V=\{v_1,\ldots,v_n\}$ and $W=\{w_1,\ldots,w_n\}$. Denote by $M$ the number of perfect matchings and
$d_{v_i}$ the degree of $v_i$.  Then \[
M\leq \prod_{i=1}^n(d(v_i)!)^{1/d(v_i)}
\]
\end{theorem}


\begin{thebibliography}{9}

\bibitem{Bregman}  L. M.  Bregman,  \emph {Certain properties of nonnegative matrices and their permanents,}   (Russian) Dokl. Akad. Nauk SSSR 211 (1973), 27--30.


\bibitem{Bourgade} P. Bourgade,  C. Hughes,  A.  Nikeghbali and M.  Yor,  \emph {The characteristic polynomial of a random unitary matrix: a probabilistic approach,}  Duke Math. J. 145 (2008), no. 1, 45--69.


\bibitem{Janson} S. Janson, \emph{ The Numbers of Spanning Trees, Hamilton Cycles and Perfect Matchings in a
Random Graph,}  Combinatorics, Probability and Computing, 3, pp 97--126.

\bibitem{JLR}  S. Janson, T. Luczak and A. Rucinski,  Random graphs, Wiley-Interscience Series in Discrete Mathematics and Optimization. Wiley-Interscience, New York, 2000.

\bibitem{Khinchin}  A. Khinchine,  \emph{\"{U}ber einen Satz der Wahrscheinlichkeitsrechnung}, Fundamenta Mathematicae 6 (1924) 9--20.

\bibitem{Kolmogorov} A. Kolmogoroff,  \emph{\"{U}ber das Gesetz des iterierten Logarithmus,}  Mathematische Annalen, 101 (1929) 126--135.

\bibitem{NV}  H. Nguyen and V. Vu,  \emph {Random matrices: law of the determinant,}  Ann. Probab. 42 (2014), no. 1, 146--167.


\bibitem{Rem}   G. Rempala  and J. Wesolowski,  \emph{Central limit theorems for random permanents with correlation structure,} J. Theoret. Probab. 15 (2002), no. 1, 63--76.




\bibitem{PZ}   E.  Paquette and O.  Zeitouni, \emph{ Extremal eigenvalue fluctuations in the GUE minor process and the law of fractional logarithm,} preprint, arXiv:1505.05627.

\bibitem{Pastur} L.  Pastur and M. Shcherbina, Eigenvalue distribution of large random matrices,  Mathematical Surveys and Monographs, 171. American Mathematical Society, Providence, RI, 2011.



\end{thebibliography}
\end{document}